\documentclass[12pt, a4paper, twoside]{article}
\usepackage{amssymb}
\usepackage[leqno]{amsmath}
\usepackage{latexsym}
\usepackage{enumerate}
\usepackage{amsmath, amssymb, amscd}
\DeclareMathOperator{\Pic}{Pic}

\DeclareMathOperator{\Cliff}{Cliff}

\begin{document}

\title{\textbf{\Large{ACM line bundles on polarized K3 surfaces}}}

\author{Kenta Watanabe \thanks{Nihon University, College of Science and Technology, 7-24-1 Narashinodai Funabashi city Chiba 274-8501 Japan, {\it E-mail address:watanabe.kenta@nihon-u.ac.jp}, Telephone numbers: 090-9777-1974} }

\date{}

\maketitle


\begin{abstract} \noindent An ACM bundle on a polarized algebraic variety is defined as a vector bundle whose intermediate cohomology vanishes. We are interested in ACM bundles of rank one with respect to a very ample line bundle on a K3 surface. In this paper, we give a necessary and sufficient condition for a non-trivial line bundle $\mathcal{O}_X(D)$ on $X$ with $|D|\neq\emptyset$ and $D^2\geq L^2-6$ to be an ACM and initialized line bundle with respect to $L$, for a given K3 surface $X$ and a very ample line bundle $L$ on $X$.\end{abstract}

\section{Introduction}

Let $X$ be a smooth projective surface. For a given very ample line bundle $L$ on $X$, we call a vector bundle $E$ on $X$ an {\it{Arithmetically Cohen-Macaulay}} ({\it{ACM}} for short) bundle with respect to $L$ if $H^1(X,E\otimes L^{\otimes l})=0$, for any integer $l\in\mathbb{Z}$. Previously, many people have investigated indecomposable ACM bundles of higher rank on several types of smooth polarized surfaces. For example, Kn${\rm{\ddot{o}}}$rrer [Kn] has proved that if $X$ is a quadric in $\mathbb{P}^3$, there are finitely many indecomposable ACM bundles of rank $r\geq2$ on $X$ with respect to the invertible sheaf defined by a hyperplane section of $X$. If $X$ is a cubic surface in $\mathbb{P}^3$, Casanellas and Hartshorne [C-H] have constructed an $n^2+1$-dimensional family of indecomposable ACM bundles of rank $n$ on $X$ with Chern classes $c_1=\mathcal{O}_X(n)$ and $c_2=\frac{1}{2}(3n^2-n)$ for $n\geq2$. 

In general, it is difficult to give a classification of indecomposable ACM bundles with respect to a given polarization. However, we can easily see that extensions of ACM vector bundles are ACM as well. Hence, we often classify ACM line bundles to construct indecomposable ACM bundles of rank $r\geq2$. For example, Joan Pons-Llopis and Fabio Tonini [P-T] have classified ACM line bundles on a DelPezzo surface $X$ with respect to the anti-canonical line bundle of $X$, and have constructed families of indecomposable ACM bundles on $X$ of rank $r\geq2$, by using extensions of ACM line bundles on $X$. On the other hand, Gianfranco Casnati [C] has classified ACM bundles of rank 2 on general determinantal quartic hypersurfaces in $\mathbb{P}^3$. This result extends our previous work [W] about the classification of ACM line bundles on quartic hypersurfaces in $\mathbb{P}^3$.

\newtheorem{thm}{Theorem}[section]

\begin{thm} {\rm{([W], Theorem 1.1)}} Let $X$ be a smooth quartic hypersurface in $\mathbb{P}^3$, and let $D$ be a nonzero effective divisor on $X$. Then the following conditions are equivalent.

\smallskip

\smallskip

\noindent {\rm{(i)}} $\mathcal{O}_X(D)$ is an ACM and initialized line bundle.

\noindent {\rm{(ii)}} For a hyperplane section $H$ of $X$, one of the following cases occurs.

\smallskip

\smallskip

{\rm{(a)}} $D^2=-2$ and $1\leq H.D\leq 3$.

{\rm{(b)}} $D^2=0$ and $3\leq H.D\leq 4$.

{\rm{(c)}} $D^2=2$ and $H.D=5$.

{\rm{(d)}} $D^2=4,\;H.D=6$ and $|D-H|=|2H-D|=\emptyset.$\end{thm}

\noindent We call the pair $(X,L)$ consisting of an algebraic K3 surface $X$ and a very ample line bundle $L$ on $X$ a polarized K3 surface, and call the sectional genus of $L$ the genus of it. We often say that a line bundle $M$ on a K3 surface $X$ is {\it{initialized}} with respect to a given polarization $L$ on $X$ if $H^0(X,M)\neq\emptyset$ and $H^0(X,M\otimes L^{\vee})=0$, and call an ACM line bundle on $X$ with respect to $L$ an ACM line bundle on $(X,L)$. In Theorem 1.1, we gave a numerical characterization of ACM and initialized line bundles on a polarized K3 surface consisting of a smooth quartic hypersurface in $\mathbb{P}^3$ and a hyperplane section of it. In this paper, we give the following result on ACM and initialized line bundles on a polarized K3 surface $(X,L)$ of any genus $g\geq 3$ as a generalization of Theorem 1.1.

\begin{thm} Let $X$ be a K3 surface, and let $L$ be a very ample line bundle on $X$. Let $D$ be a nonzero effective divisor on $X$ with $D^2\geq L^2-6$. Then the following conditions are equivalent.

\smallskip

\smallskip

\noindent {\rm{(i)}} $\mathcal{O}_X(D)$ is an ACM and initialized line bundle with respect to $L$.

\noindent {\rm{(ii)}} For $H\in |L|$, one of the following cases occurs.

\smallskip

\smallskip

{\rm{(a)}} $D^2 = H^2-6$ and $H^2-3\leq H.D\leq H^2-1.$

{\rm{(b)}} $D^2 = H^2-4$ and $H^2-1\leq H.D\leq H^2.$

{\rm{(c)}} $D^2 = H^2-2$ and $H.D=H^2+1.$

{\rm{(d)}} $D^2\geq H^2,\;D^2=2H.D-H^2-4,\;|D-H|=\emptyset$ and $h^1(\mathcal{O}_X(2H-D))=0.$

\end{thm}

\noindent In Theorem 1.2, since $H$ is not hyperelliptic, $H^2\geq4$. Moreover, the condition $h^1(\mathcal{O}_X(2H-D))=0$ as in Theorem 1.2 (ii) (d) implies that $$h^0(\mathcal{O}_X(2H-D))=\displaystyle\frac{3}{2}H^2-H.D.$$
In particular, we note that if $H^2=4$, it is equivalent to the condition $|2H-D|=\emptyset$ as in Theorem 1.1 (ii) (d).

Our plan of this paper is as follows. In Section 2, we recall some basic results about line bundles and linear systems on K3 surfaces. In Section 3, we give a numerical characterization of ACM line bundles on polarized K3 surfaces, and prove our main theorem. In Section 4, we give an example of an ACM line bundle on certain polarized K3 surfaces of Picard number 2.

\smallskip

\smallskip

\noindent {\bf{Notation and Conventions.}} We work over the complex number field $\mathbb{C}$. A surface is a smooth projective surface. Let $X$ be a surface. We denote by $\Pic(X)$ the Picard group of $X$. For a divisor $D$ on $X$, we will denote by $|D|$ the linear system defined by $D$. If two divisors $D_1$ and $D_2$ on $X$ are linearly equivalent, then we will write $D_1\sim D_2$. We call a regular surface a K3 surface if the canonical line bundle of it is trivial.

\section{Linear systems and line bundles on K3 surfaces}

In this section, we recall some basic results about ample line bundles and linear systems on K3 surfaces. First of all, we remark some facts about numerically connected divisors on a surface.

\newtheorem{df}{Definition}[section]

\begin{df} A non-zero effective divisor $D$ on a surface is called {\rm{$m$-connected}} if $D_1.D_2\geq m$, for each effective decomposition $D=D_1+D_2$.\end{df}

\noindent If a non-zero effective divisor $D$ on a surface is 1-connected, then $h^0(\mathcal{O}_{D})=1$ (cf. [B-P-W], Corollary 12.3). Hence, we can easily see that, for a 1-connected divisor $D$ on a K3 surface, we get $h^1(\mathcal{O}_X(D))=0$. Next, we recall a result about the classification of base point free divisors on K3 surfaces.

\newtheorem{prop}{Proposition}[section]

\begin{prop}{\rm{([SD], 2.7)}} Let $L$ be a numerically effective line bundle on a K3 surface $X$. Then $|L|$ is not base point free if and only if there exist an elliptic curve $F$, a smooth rational curve $\Gamma$ and an integer $k\geq2$ such that $F.\Gamma=1$ and $L\cong\mathcal{O}_X(kF+\Gamma)$. \end{prop}

\begin{prop}{\rm{([SD], Proposition 2.6)}} Let $L$ be a line bundle on a K3 surface $X$ such that $|L|\neq\emptyset$. Assume that $|L|$ has no fixed component. Then one of the following cases occurs.

\smallskip

\smallskip

{\rm{(i)}} $L^2>0$ and the general member of $|L|$ is a smooth irreducible curve of genus $\frac{1}{2}L^2+1$.

{\rm{(ii)}} $L^2=0$ and $L\cong\mathcal{O}_X(kF)$, where $k\geq1$ is an integer and $F$ is a smooth curve of genus one. In this case, $h^1(L)=k-1$. \end{prop} 

\noindent It is well known that, for an irreducible curve $C$ on a K3 surface such that $C^2>0$, $|C|$ is base point free ([SD], Theorem 3.1). Hence, by Proposition 2.2, the following proposition follows.

\begin{prop}{\rm{([SD], Corollary 3.2)}} Let $L$ be a line bundle such that $|L|\neq\emptyset$ on a K3 surface. Then $|L|$ has no base point outside its fixed component.\end{prop}

\noindent {\bf{Remark 2.1}} It is well known that, if $X$ is a K3 surface, then the self-intersection of any divisor on $X$ is an even integer. In particular, any rational curve $\Delta$ on $X$ satisfies $\Delta^2=-2$. Moreover, the self-intersection of the fixed component of any non-zero effective divisor is a negative integer.

$\;$

\noindent At the end of this section, we recall some classical results about very ample line bundles on K3 surfaces. It is well known that if an ample linear system on a K3 surface is not very ample, then it is hyperelliptic ([SD]). Hence, by the characterization of hyperelliptic linear systems on K3 surfaces, we have the following assertion.

\begin{prop} {\rm{(cf. [M-M], and [SD], Theorem 5.2)}} Let $L$ be a numerically effective line bundle with $L^2\geq4$ on a K3 surface $X$. Then $L$ is very ample if and only if the following conditions are satisfied.

\smallskip

\smallskip

{\rm{(i)}} There is no irreducible curve $E$ such that $E^2=0$ and $E.L=1$ or 2.

{\rm{(ii)}} There is no irreducible curve $E$ such that $E^2=2$ and $L\cong\mathcal{O}_X(2E)$.

{\rm{(iii)}} There is no irreducible curve $E$ such that $E^2=-2$ and $E.L=0$. \end{prop}

\noindent Note that, by Proposition 2.1 and Proposition 2.4, if $L$ is a very ample line bundle, then $|L|$ is base point free. Hence, the general member of it is a smooth irreducible curve. Moreover, by Proposition 2.2, we have the following fact.

$\;$

\noindent {\bf{Corollary 2.1}} {\it{Let $X$ be a K3 surface, and $L$ be a very ample line bundle on $X$. Then the following statements hold.

\smallskip

\smallskip

{\rm{(i)}} Let $D$ be a nonzero effective divisor on X with $D^2\geq0$. Then $L.D\geq 3$. In particular, if $D^2=0$ and $L.D=3$, then $|D|$ is an elliptic pencil.

\smallskip

\smallskip

{\rm{(ii)}} There is no effective divisor $D$ on $X$ with $D^2=2$ and $L=\mathcal{O}_X(2D)$.}}

$\;$

{\it{Proof.}} Note that, since $L$ is very ample, we have $L^2\geq4$. 

\smallskip

(i) Let $\Delta$ be the fixed component of $|D|$. Since $D^2\geq0$, we have $D-\Delta\neq 0$. If $(D-\Delta)^2=0$, then, by Proposition 2.2 (ii), there exists an elliptic curve $F$ with $D-\Delta \sim kF\;(k\geq1)$. In this case, by Proposition 2.4, we have $L.(D-\Delta)\geq 3k\geq3$. If $(D-\Delta)^2\geq 2$, by the Hodge index theorem, we have $(L.(D-\Delta))^2\geq (L^2)(D-\Delta)^2\geq8$. Hence, we have $L.D\geq 3$. In particular, if $L.D=3$, then $\Delta$ is empty, by the very ampleness of $L$ and the above observation. Hence, in this case, if $D^2=0$, then $|D|$ is an elliptic pencil. 

(ii) Assume that $D$ is an effective divisor on $X$ with $D^2=2$ and $L=\mathcal{O}_X(2D)$. We show that $|D|$ is base point free. Let $\Delta$ be the fixed component of $|D|$ and assume that it is not empty. Since $L.D=4$, by the assertion of (i), we have $L.(D-\Delta)=3$. Hence, $\Delta$ is irreducible. Since $L^2=8$, by the Hodge index theorem, we have $(D-\Delta)^2=0$. This implies that $|D-\Delta|$ is an elliptic pencil on $X$ with $(D-\Delta).\Delta=2$. However, this means that $\Delta$ is not the fixed component of $|D|$. Therefore, $|D|$ is base point free. By Proposition 2.2 (i) and Proposition 2.4, this is a contradiction. $\hfill\square$

$\;$

\noindent {\bf{Remark 2.2}}. If $(X,L)$ is a polarized K3 surface, for a non-zero effective divisor $D$ on $X$, $\mathcal{O}_X(D)$ is ACM if and only if the same is true for $\mathcal{O}_X(-D)$. 

\section{Proof of Theorem 1.2}

In this section, we give a proof of Theorem 1.2. First of all, in order to prove our main theorem, we prove the following lemmas.

\newtheorem{lem}{Lemma}[section]

\begin{lem} Let $X$ be a K3 surface, let $L$ be a very ample line bundle on $X$, and let $D$ be a nonzero effective divisor on $X$. Moreover, let $m\in\mathbb{N}$. If $L.D\leq mL^2-1$ and, for any $k\in\mathbb{Z}$ with $0\leq k\leq m$,  $h^1(\mathcal{O}_X(D)\otimes L^{\vee\otimes k})=0$, then $\mathcal{O}_X(D)$ is an ACM line bundle.\end{lem}

{\it{Proof}}. Let $m$ be a positive integer satisfying the assumption. Let $n\in\mathbb{N}$ and let $H\in |L|$ be a smooth irreducible curve. First of all, we have
$$h^1(\mathcal{O}_H(nH+D))=h^0(\mathcal{O}_H(-D-(n-1)H))=0.$$
\noindent The above vanishing and the cohomology of the exact sequence
$$0\rightarrow\mathcal{O}_X(D+(n-1)H)\rightarrow\mathcal{O}_X(D+nH)\rightarrow\mathcal{O}_H(D+nH)\rightarrow0$$
\noindent imply 
$$h^1(\mathcal{O}_X(D+nH))\leq h^1(\mathcal{O}_X(D+(n-1)H)).$$
\noindent Since $h^1(\mathcal{O}_X(D))=0$, it follows that $h^1(\mathcal{O}_X(D+nH))=0$,
by descending induction on $n$.

On the other hand, since $L.D\leq mL^2-1$, if $n\geq m$, then we have
$$h^1(\mathcal{O}_H((n+1)H-D))=h^0(\mathcal{O}_H(D-nH))=0.$$
\noindent By the exact sequence
$$0\rightarrow\mathcal{O}_X(nH-D)\rightarrow\mathcal{O}_X((n+1)H-D)\rightarrow\mathcal{O}_H((n+1)H-D)\rightarrow0,$$
\noindent and induction on $n$, we have
$$h^1(\mathcal{O}_X(D-(n+1)H))=h^1(\mathcal{O}_X((n+1)H-D))=0\;(n\geq m).$$
\noindent Hence, $\mathcal{O}_X(D)$ is ACM.$\hfill\square$

\begin{lem} Let $X$ be a K3 surface, and let $D$ be an effective divisor on $X$ which is not linearly equivalent to 0. Let $\Delta$ be the fixed component of $|D|$. If $h^1(\mathcal{O}_X(D-\Delta))=0$ and $D^2=(D-\Delta)^2$, then $h^1(\mathcal{O}_X(D))=0$.\end{lem}

{\it{Proof}}. Since $\Delta$ is the fixed component of $|D|$, it follows that $\Delta^2<0$. If $\Delta\sim D$, then we have the contradiction
$$0=(D-\Delta)^2=D^2=\Delta^2<0.$$ 
Hence, the movable part of $|D|$ is not empty. Since $h^1(\mathcal{O}_X(D-\Delta))=0$, we have
$$h^0(\mathcal{O}_X(D))=h^0(\mathcal{O}_X(D-\Delta))=\chi(\mathcal{O}_X(D-\Delta)).$$
\noindent On the other hand, since $D^2=(D-\Delta)^2$, we have $\chi(\mathcal{O}_X(D))=\chi(\mathcal{O}_X(D-\Delta))$. Hence, we have $h^1(\mathcal{O}_X(D))=0$. $\hfill\square$

$\;$

{\it{Proof of Theorem 1.2}}. Let $X$ be a K3 surface, $L$ be a very ample line bundle, and let $D$ be a non-zero effective divisor on $X$ with $D^2\geq L^2-6$. Let $H\in|L|$ be a smooth curve. If $H^2=4$, then the statement is already proved in Theorem 1.1. Hence, from now on we assume that $H^2\geq6$. Since the proof is so long and complex, we divide it into several cases.

\begin{prop} Assume that $D^2=H^2-6$. Then $\mathcal{O}_X(D)$ is ACM and initialized with respect to $L$ if and only if $H^2-3\leq H.D\leq H^2-1$. \end{prop}

{\it{Proof}}. Assume that $\mathcal{O}_X(D)$ is ACM and initialized. If $|H-D|=\emptyset$, then $\chi(\mathcal{O}_X(H-D))=0$. Hence, we have $H.D=H^2-1$. We consider the case where $|H-D|\neq\emptyset$. Since $H.(H-D)\geq1$, we have 
$$H.D\leq H^2-1.\leqno(1)$$ 
If $H^2=6$, then $D^2=0$. Hence, we have $H.D\geq 3=H^2-3$, by Corollary 2.1 (i). Assume that $H^2=8$. Then we have $D^2=2$. By the Hodge index theorem, we have $H.D\geq4$. If $H.D=4$, we have $(H-2D)^2=0$ and $H.(H-2D)=0$. By the ampleness of $H$, we have $H\sim 2D$. By Corollary 2.1 (ii), this is a contradiction. Therefore, we have $H.D\geq 5=H^2-3$. If $H^2\geq 10$, by the Hodge index theorem, we have
$$(H^2-4)^2<H^2(H^2-6)=H^2D^2\leq (H.D)^2.$$
Thus, $H^2-3\leq H.D$. Moreover, by the inequality (1), we have
$$H^2-3\leq H.D\leq H^2-1.\leqno(2)$$

$\;$

Conversely, we assume that the inequality (2) holds. Since $H.(D-H)\leq -1$, we have $|D-H|=\emptyset$. Hence, $\mathcal{O}_X(D)$ is initialized. We show that it is ACM with respect to $H$. Since $H.D\leq H^2-1$, by Lemma 3.1, it is sufficient to show that
$$h^1(\mathcal{O}_X(D))=h^1(\mathcal{O}_X(H-D))=0.$$
First of all, we show that $h^1(\mathcal{O}_X(H-D))=0$. If $H.D=H^2-3$, then we have
$$(H-D)^2=0\text{ and } H.(H-D)=3.$$
By Corollary 2.1 (i), $|H-D|$ is an elliptic pencil. Therefore, $h^1(\mathcal{O}_X(H-D))=0$. If $H.D=H^2-2$, we have $(H-D)^2=-2$. Since $H.(H-D)=2$, by Corollary 2.1 (i), the movable part of $|H-D|$ is empty. This implies that $h^0(\mathcal{O}_X(H-D))=1$, and hence, we have $h^1(\mathcal{O}_X(H-D))=0$. Assume that $H.D=H^2-1$. Since $(H-D)^2=-4$, we have $|H-D|=\emptyset$. In fact, since $H.(H-D)=1$, if $|H-D|\neq\emptyset$, the member of it is a $(-2)$-curve. This is a contradiction. Since $|D-H|=\emptyset$, we have $h^1(\mathcal{O}_X(H-D))=0$.

Next we show that $h^1(\mathcal{O}_X(D))=0$. Assume that $|D|$ is base point free. If $H^2\geq8$, then $D^2\geq2$. Hence, by the theorem of Bertini, we have $h^1(\mathcal{O}_X(D))=0$. Assume that $H^2=6$. Since $D^2=0$, by Proposition 2.2 (ii), there exist an elliptic curve $F$ on $X$ and a positive integer $k$ such that $D\sim kF$. By Corollary 2.1 (i), we have $H.D\geq 3k$. Since $H.D\leq H^2-1=5$, we have $k=1$. Hence, $h^1(\mathcal{O}_X(D))=0$. Assume that $|D|$ is not base point free. Let $\Delta$ be the fixed component of $|D|$ and let $D^{'}$ be the movable part of $|D|$. Then since $H.D\leq H^2-1$, by the ampleness of $H$, we have 
$$H.D^{'}\leq H^2-2.$$

We consider the case where ${D^{'}}^2=0$. Then we have $H^2=6\text{ or }8$. Indeed, by Proposition 2.2 (ii), there exist an elliptic curve $F$ and a positive integer $k$ such that $D^{'}\sim kF$. By Corollary 2.1 (i), $H.F\geq3$. Hence, we have $k\leq\displaystyle\frac{1}{3}(H^2-2)$. Since $h^1(\mathcal{O}_X(D^{'}))=k-1$, we have
$$\chi(\mathcal{O}_X(D^{'}))+k-1=h^0(\mathcal{O}_X(D^{'}))=h^0(\mathcal{O}_X(D))\geq\chi(\mathcal{O}_X(D)).$$
This implies that ${D^{'}}^2\geq\displaystyle\frac{1}{3}(H^2-8)$. If $H^2\geq 10$, this is a contradiction.

Assume that $H^2=6$. Since $D^2=0$ and $H.D^{'}\leq 4$, by Corollary 2.1 (i), we have $k=1$. Hence, $|D^{'}|$ is an elliptic pencil. By Lemma 3.2, we have $h^1(\mathcal{O}_X(D))=0$.

Assume that $H^2=8$. Since $H.D^{'}\leq 6$, we have $k=2$. Indeed, since $D^2=2$, we have $h^0(\mathcal{O}_X(D))\geq 3$. If $k=1$, we have the contradiction $$h^0(\mathcal{O}_X(D))=h^0(\mathcal{O}_X(D^{'}))=2.$$
Since $H.D^{'}=6$, we have $H.\Delta=1$. This implies that $\Delta$ is a $(-2)$-curve. Since $D^2=2$, we have $F.\Delta=1$. Hence, $h^1(\mathcal{O}_X(D))=0$.

We consider the case where ${D^{'}}^2>0$. Since $h^1(\mathcal{O}_X(D^{'}))=0$, we have
$$\chi(\mathcal{O}_X(D^{'}))=h^0(\mathcal{O}_X(D^{'}))=h^0(\mathcal{O}_X(D))\geq\chi(\mathcal{O}_X(D)).$$
Hence, we have ${D^{'}}^2\geq H^2-6$. Assume that ${D^{'}}^2\geq H^2-4$. Since $H^2\geq6$, we have ${D^{'}}^2\geq2$. It follows from the Hodge index theorem that
$$(H^2-3)^2<H^2(H^2-4)\leq H^2{D^{'}}^2\leq (H.D^{'})^2.$$
Thus, we have $H^2-2\leq H.D^{'}$. Since $H.D^{'}\leq H^2-2$, we have $H.D^{'}=H^2-2$. However, since $(H-D^{'})^2\geq0$ and $H.(H-D^{'})=2$, by Corollary 2.1 (i), this is a contradiction. Therefore, we have $D^2=H^2-6={D^{'}}^2>0$. Since $h^1(\mathcal{O}_X(D^{'}))=0$, by Lemma 3.2, we have $h^1(\mathcal{O}_X(D))=0$. $\hfill\square$

\begin{prop} Assume that $D^2=H^2-4$. Then $\mathcal{O}_X(D)$ is ACM and initialized with respect to $L$ if and only if $H.D=H^2-1$ or $H^2$.\end{prop}

{\it{Proof}}. Assume that $\mathcal{O}_X(D)$ is ACM and initialized. If $|H-D|=\emptyset$, by the assumption, we have $\chi(\mathcal{O}_X(H-D))=0$. Hence, we have $H.D=H^2$. Assume that $|H-D|\neq\emptyset$. By the ampleness of $H$, we have $H.(H-D)\geq1$, that is, $H.D\leq H^2-1$. On the other hand, since $H^2\geq6$, by the Hodge index theorem, we have
$$(H^2-3)^2<H^2(H^2-4)=H^2D^2\leq(H.D)^2.$$
Hence, we have $H.D=H^2-1$. Indeed, since $H^2-3<H.D$, we have 
$$H.D=H^2-2\text{ or }H^2-1.$$
If $H.D=H^2-2$, we have $(H-D)^2=0$ and $H.(H-D)=2$. However, this contradicts Corollary 2.1 (i).

Conversely, we assume that $H.D=H^2-1$ or $H^2$. \noindent Since $H.(D-H)=-1$ or 0, we have $|D-H|=\emptyset$. Hence, $\mathcal{O}_X(D)$ is initialized. We show that $\mathcal{O}_X(D)$ is ACM with respect to $H$.

Assume that $H.D=H^2-1$. By Lemma 3.1, it is sufficient to show that 
$$h^1(\mathcal{O}_X(D))=h^1(\mathcal{O}_X(H-D))=0.$$
First of all, since $(H-D)^2=-2$ and $H.(H-D)=1$, the member of $|H-D|$ is irreducible. Therefore, we have $h^1(\mathcal{O}_X(H-D))=0$.
\noindent In order to show that $h^1(\mathcal{O}_X(D))=0$, we show that $|D|$ is base point free. Assume that $|D|$ is not base point free. Let $\Delta$ be the fixed component of $|D|$, and let $D^{'}$ be the movable part of $|D|$. Then we note that, by the ampleness of $H$, we have 
$$H.D^{'}\leq H^2-2.\leqno(3)$$
Assume that ${D^{'}}^2=0.$ Then, by Proposition 2.2 (ii), there exist an elliptic curve $F$ and an integer $k\geq1$ such that $D^{'}\sim kF$. By the inequality (3) and Proposition 2.4 (i), we have $k\leq\frac{1}{3}(H^2-2)$. Since, by Proposition 2.2 (ii), $h^1(\mathcal{O}_X(D^{'}))=k-1$, we have
$$\chi(\mathcal{O}_X(D^{'}))\geq\chi(\mathcal{O}_X(D))-\frac{1}{3}(H^2-5).$$
\noindent Since $H^2\geq6$, we have the contradiction
$${D^{'}}^2\geq \frac{1}{3}(H^2-2)>0.$$
\noindent Since ${D^{'}}^2>0$, we have $h^1(\mathcal{O}_X(D^{'}))=0$. Hence, by comparing $\chi(\mathcal{O}_X(D^{'}))$ and $\chi(\mathcal{O}_X(D))$, we have ${D^{'}}^2\geq D^2=H^2-4$. Since $H^2\geq6$, by the Hodge index theorem, we have
$$(H^2-3)^2<H^2(H^2-4)\leq (H^2)({D^{'}}^2)\leq (H.D^{'})^2.$$
Thus, $H^2-3<H.D^{'}$.
\noindent By the inequality (3), we have $H.D^{'}=H^2-2$. Hence, we also have $(H-D^{'})^2\geq0$ and $H.(H-D^{'})=2$. However, this contradicts Corollary 2.1 (i). Therefore, $|D|$ is base point free. Since $H^2\geq6$, we have $D^2>0$. This implies that $h^1(\mathcal{O}_X(D))=0$.

Assume that $H.D=H^2$. Since $H.(D-H)=0$, by the ampleness of $H$, we have $|H-D|=\emptyset$. Since $D^2=H^2-4$, we have $(H-D)^2=-4$. Hence, we have 
$$h^1(\mathcal{O}_X(H-D))=-\chi(\mathcal{O}_X(H-D))=0.$$
\noindent By Lemma 3.1, it is sufficient to show that 
$$h^1(\mathcal{O}_X(D))=h^1(\mathcal{O}_X(2H-D))=0.$$ 

First of all, we consider the case where $|D|$ is base point free. Since $D^2>0$, we have $h^1(\mathcal{O}_X(D))=0$. In order to show that $h^1(\mathcal{O}_X(2H-D))=0$, we show that $|2H-D|$ is base point free. Assume that it is not base point free, and let $\Delta$ be the fixed component of it. Then, since $(2H-D)^2=D^2>0$ and $H.(2H-D)=H^2>0$, the movable part of $|2H-D|$ is not empty. Hence, we take a nonzero divisor $D^{'}\in|2H-D-\Delta|$. Note that since $H.(2H-D)=H^2$, by the ampleness of $H$, we have 
$$H.D^{'}\leq H^2-1.\leqno(4)$$
If ${D^{'}}^2=0$, then there exist an elliptic curve $F$ and a positive integer $k$ such that $D^{'}\sim kF$. By the inequality (4) and Proposition 2.4 (i), we have $k\leq \displaystyle\frac{1}{3}(H^2-1)$. By the same reason as above, we have
$$\chi(\mathcal{O}_X(D^{'}))\geq \chi(\mathcal{O}_X(2H-D))-\displaystyle\frac{1}{3}(H^2-4).$$
\noindent Since $H^2\geq6$, we have the contradiction ${D^{'}}^2\geq \displaystyle\frac{1}{3}(H^2-4)>0.$
Hence, we have ${D^{'}}^2>0$. By the same way as above, we have 
$${D^{'}}^2\geq (2H-D)^2=H^2-4.$$
Since $H^2\geq6$, by the Hodge index theorem, we have
$$(H^2-3)^2<(H^2)(H^2-4)\leq (H^2)({D^{'}}^2)\leq (H.D^{'})^2.$$
\noindent Hence, by the inequality (4), we have $H.D^{'}=H^2-2$ or $H^2-1$. If $H.D^{'}=H^2-2$, then we have $(H-D^{'})^2\geq 0$ and $H.(H-D^{'})=2$. However, by Corollary 2.1 (i), this is a contradiction. Therefore, we have 
$$H.D^{'}=H^2-1.\leqno(5)$$
 Note that, since $H.\Delta=1$, $\Delta$ is a $(-2)$-curve. Then we have ${D^{'}}^2=H^2-4.$ In fact, if ${D^{'}}^2>H^2-4$, by the equality (5), we have $(H-D^{'})^2>-2$ and $H.(H-D^{'})=1$. This contradicts Corollary 2.1 (i). Therefore, we have $(H-D^{'})^2=-2$. Since $H.(H-D^{'})=1$, the member of $|H-D^{'}|$ is a $(-2)$-curve. Since 
$$(D^{'}+\Delta)^2=(2H-D)^2=H^2-4,$$
\noindent we have $D^{'}.\Delta=1$. Hence, 
$$D^{'}.(2H-D)=D^{'}.(D^{'}+\Delta)=H^2-3.$$
By the equality (5), we have $D^{'}.D=H^2+1$. Thus, we have $D.(H-D^{'})=-1<0$. This contradicts the assumption that $|D|$ is base point free. Hence, $|2H-D|$ is base point free. Since $(2H-D)^2=D^2>0$, we have $h^1(\mathcal{O}_X(2H-D))=0$.

We consider the case where $|D|$ is not base point free. Let $\Delta$ be the fixed component of $|D|$, and let $D^{'}$ be the movable part of $|D|$. Note that, since $H.D=H^2$, we have 
$$H.D^{'}\leq H^2-1.\leqno(6)$$ 
Since $H^2\geq6$, by the same reason as above, we have ${D^{'}}^2>0$ and hence, we have $h^1(\mathcal{O}_X(D^{'}))=0$. Since $\chi(\mathcal{O}_X(D^{'}))\geq \chi(\mathcal{O}_X(D))$, we have ${D^{'}}^2\geq D^2=H^2-4$. Since $H^2\geq6$, by the Hodge index theorem, we have
$$(H^2-3)^2<(H^2)(H^2-4)\leq(H^2)({D^{'}}^{2})\leq (H.D^{'})^2.$$
Thus, $H^2-3<H.D^{'}.$
\noindent Therefore, by the inequality (6), we have $H.D^{'}=H^2-2$ or $H^2-1$. By Corollary 2.1 (i), if $H.D^{'}=H^2-2$, we have the contradiction $(H-D^{'})^2\geq 0$ and $H.(H-D^{'})=2$. Hence, we have 
$$H.D^{'}=H^2-1.\leqno(7)$$ 
Since $H.\Delta=1$, $\Delta$ is a $(-2)$-curve. Hence, $D$ is a 1-connected divisor. Indeed, since $D^2=H^2-4$, we have $2D^{'}.\Delta=H^2-2-{D^{'}}^2$. Since $D^{'}.\Delta\geq0$, we have ${D^{'}}^2=H^2-2$ or $H^2-4$. If ${D^{'}}^2=H^2-2$, by the equality (7), we have $(H-D^{'})^2=0$ and $H.(H-D^{'})=1$. However, by Corollary 2.1 (i), this is a contradiction. Hence, we have ${D^{'}}^2=H^2-4$ and hence, we have 
$$D^{'}.\Delta=1.\leqno(8)$$
Since we may assume that $D^{'}$ is irreducible, $|D|$ contains a 1-connected divisor.
Therefore, we have $h^1(\mathcal{O}_X(D))=0$. 

Next, we show that $h^1(\mathcal{O}_X(2H-D))=0$. Let $D^{'}$ and $\Delta$ be as above. Note that 
$$(H-\Delta)+(H-D^{'})=2H-D.$$ 
By the equality (7) and (8), we have $(H-\Delta).(H-D^{'})=1$. Moreover, since $H.(H-D^{'})=1$ and $(H-D^{'})^2=-2$, the member of $|H-D^{'}|$ is a $(-2)$-curve. Here, in order to show that $|2H-D|$ contains a 1-connected divisor, we show that $|H-\Delta|$ is base point free. We assume that it is not base point free and let $\Delta^{'}$ be the fixed component of it. Since $(H-\Delta)^2>0$, the movable part of $|H-\Delta|$ is not empty. Hence, we take a nonzero divisor $D^{''}\in |H-\Delta-\Delta^{'}|$. Note that, since $H.\Delta=1$, we have 
$$H.D^{''}\leq H^2-2.\leqno(9)$$
Then we have ${D^{''}}^2>0$. In fact, if ${D^{''}}^2=0$, then there exist a positive integer $k$ and an elliptic curve $F$ on $X$ such that $D^{''}\sim kF$. By the inequality (9) and Corollary 2.1 (i), we have
$$k\leq\displaystyle\frac{1}{3}(H^2-2).$$
Since 
$$\chi(\mathcal{O}_X(D^{''}))+k-1=h^0(\mathcal{O}_X(D^{''}))=h^0(\mathcal{O}_X(H-\Delta))\geq\chi(\mathcal{O}_X(H-\Delta)),$$
and $H^2\geq6$, we have the contradiction
$${D^{''}}^2\geq\displaystyle\frac{1}{3}(H^2-2)>0.$$
 Hence, we have $h^1(\mathcal{O}_X(D^{''}))=0$. By comparing $\chi(\mathcal{O}_X(D^{''}))$ and $\chi(\mathcal{O}_X(H-\Delta))$, we have ${D^{''}}^2\geq (H-\Delta)^2=H^2-4$. It follows from the Hodge index theorem, that
$$(H^2-3)^2<(H^2)(H^2-4)\leq(H^2)({D^{''}}^2)\leq (H.D^{''})^2.$$
\noindent Hence, by the inequality (9), we have $H.D^{''}=H^2-2$ and hence, we have $H.(H-D^{''})=2$. Since $(H-D^{''})^2\geq0$, this contradicts Corollary 2.1 (i). Hence, $|H-\Delta|$ is base point free. Since $(H-\Delta)^2>0$, the general member of it is irreducible. Therefore, $|2H-D|$ contains a 1-connected divisor. Hence, we have $h^1(\mathcal{O}_X(2H-D))=0$. $\hfill\square$

\begin{prop} Assume that $D^2=H^2-2$. Then $\mathcal{O}_X(D)$ is ACM and initialized with respect to $L$ if and only if $H.D=H^2+1$.\end{prop}

{\it{Proof}}. Assume that $\mathcal{O}_X(D)$ is ACM and initialized. Then we have $|H-D|=\emptyset$. In fact, if $|H-D|\neq\emptyset$, by the ampleness of $H$, we have $H.(H-D)\geq1$. Therefore, we have $D.H\leq H^2-1$. Since $D^2=H^2-2$, we have $(H-D)^2\geq0$. By Corollary 2.1 (i), we have $H.(H-D)\geq3$, that is, 
$$H.D\leq H^2-3.\leqno(10)$$
It follows from the Hodge index theorem that
$$(H^2-2)^2<H^2(H^2-2)=H^2D^2\leq(H.D)^2.$$
However, by the inequality (10), this is a contradiction. Therefore, by the assumption, we have $\chi(\mathcal{O}_X(H-D))=0$, and hence, $H.D=H^2+1$.

Conversely, we assume that $H.D=H^2+1$. We show that $\mathcal{O}_X(D)$ is ACM. By Lemma 3.1, it is sufficient to show that 
$$h^1(\mathcal{O}_X(D))=h^1(\mathcal{O}_X(D-H))=h^1(\mathcal{O}_X(D-2H))=0.$$
First of all, we show that $|D|$ is base point free. Assume that it is not base point free. Let $\Delta$ be the fixed component of $|D|$, and let $D^{'}$ be the movable part of $|D|$. Then we note that 
$$H.D^{'}\leq H^2\leqno(11)$$
If ${D^{'}}^2=0$, then there exist an elliptic curve $F$ and a positive integer $k$ such that $D^{'}\sim kF$. By the inequality (11) and Proposition 2.4 (i), we have
$$\chi(\mathcal{O}_X(D^{'}))\geq\chi(\mathcal{O}_X(D))-\displaystyle\frac{H^2}{3}+1.$$
Hence, we have the contradiction ${D^{'}}^2\geq\displaystyle\frac{H^2}{3}>0$. Therefore, we have ${D^{'}}^2>0$. This implies that $\chi(\mathcal{O}_X(D^{'}))\geq\chi(\mathcal{O}_X(D))$. Hence, we have ${D^{'}}^2\geq D^2=H^2-2$. By the Hodge index theorem, we have
$$(H^2-2)^2<H^2(H^2-2)\leq (H^2)({D^{'}}^2)\leq (H.D^{'})^2.$$
\noindent Hence, by the inequality (11), we have $H.D^{'}=H^2-1$ or $H^2$. If $H.D^{'}=H^2-1$, then $(H-D^{'})^2\geq 0$ and $H.(H-D^{'})=1$. This contradicts Corollary 2.1 (i). Hence, we have $H.D^{'}=H^2$. Note that, since $H.\Delta=1$, $\Delta$ is a $(-2)$-curve. Since $(H-D^{'})^2\geq -2$ and $H.(H-D^{'})=0$, we have $H\sim D^{'}$. However, since $D^2=H^2-2$, this is a contradiction. Therefore, $|D|$ is base point free. Since $H^2\geq6$ and $D^2=H^2-2$, we have $h^1(\mathcal{O}_X(D))=0$. Since $(H-D)^2=-4$, we have $\chi(\mathcal{O}_X(H-D))=0$. Since $|D|$ is base point free and $D.(D-H)=-3$, we have $|D-H|=\emptyset$. Hence, $\mathcal{O}_X(D)$ is initialized. Moreover, since $H.(H-D)=-1$, we have $|H-D|=\emptyset$. Hence, we have $h^1(\mathcal{O}_X(D-H))=0$.

We show that 
$$h^1(\mathcal{O}_X(2H-D))=0.\leqno(12)$$
If $|2H-D|$ is base point free, the equality (12) is satisfied. Indeed, if $H^2\geq8$, we have $(2H-D)^2>0$. Hence, by the theorem of Bertini, we have the the equality (12). If $H^2=6$, then $(2H-D)^2=0$. Hence, by Proposition 2.2 (ii), there exist an elliptic curve $F$ and an integer $k\geq1$ such that $2H-D\sim kF$. In this case, since $H.(2H-D)=5$, by Proposition 2.4 (i), we have $k=1$. Hence, the equality (12) holds. Therefore, we assume that it is not base point free. Let $\Delta$ be the fixed component of $|2H-D|$, and let $D^{'}$ be the movable part of $|2H-D|$. Note that, since $(2H-D)^2\geq0$, we have $D^{'}\neq0$.

First of all, we consider the case where $H^2=6$. Since $H.(2H-D)=5$, we have $H.D^{'}\leq 4$. Then we have ${D^{'}}^2=0$. In fact, if ${D^{'}}^2>0$, by the Hodge index theorem, we have
$$(H.D^{'})^2\geq (H^2)({D^{'}})^2\geq 12.$$
\noindent Thus, we have $H.D^{'}=4$. Since $(H-D^{'})^2\geq0$ and $H.(H-D^{'})=2$, by Corollary 2.1 (i), we have a contradiction. Hence, there exist an elliptic curve $F$ and an integer $k\geq1$ such that $D^{'}\sim kF$. Hence, by Corollary 2.1 (i), we have $3k\leq H.D^{'}\leq 4$, and hence, we have $k=1$. Since $h^1(\mathcal{O}_X(D^{'}))=0$ and $(2H-D)^2={D^{'}}^2$, by Lemma 3.2, the equality (12) holds.

Next, we consider the case where $H^2=8$. Since $H.(2H-D)=7$, we have 
$$H.D^{'}\leq 6.\leqno(13)$$ 

If ${D^{'}}^2=0$, there exist an elliptic curve $F$ and an integer $k\geq1$ such that $D^{'}\sim kF$. We have $3k\leq H.D^{'}\leq 6$. Hence, we have $k=1$ or 2. If $k=1$, we have $h^1(\mathcal{O}_X(D^{'}))=0$ and hence, $\chi(\mathcal{O}_X(D^{'}))\geq \chi(\mathcal{O}_X(2H-D))$. This implies the contradiction that
$${D^{'}}^2\geq (2H-D)^2=2.$$
 If $k=2$, we have $H.\Delta=1$ and hence, $\Delta$ is a $(-2)$-curve. Since $$(D^{'}+\Delta)^2=(2H-D)^2=2,$$
 we have $D^{'}.\Delta=2$. Hence, we have the equality (12). 

Assume that ${D^{'}}^2>0$. Then we have ${D^{'}}^2=2$. In fact, if we assume that ${D^{'}}^2\geq 4$, by the Hodge index theorem, and the assumption that $H^2=8$, we have
$$32\leq(H^2)({D^{'}}^2)\leq (H.D^{'})^2.$$
 Since the inequality (13) implies that $H.D^{'}=6$, we have 
$$(H-D^{'})^2\geq 0\text{ and }H.(H-D^{'})=2.$$
By Corollary 2.1 (i), this is a contradiction. Hence, we have
$$h^1(\mathcal{O}_X(D^{'}))=0\text{ and }{D^{'}}^2=(2H-D)^2=2.$$
 By Lemma 3.2, the equality (12) holds.

Finally, we consider the case where $H^2\geq 10$. We note that, since $H.(2H-D)=H^2-1$, we have 
$$H.D^{'}\leq H^2-2.\leqno(14)$$
Assume that ${D^{'}}^2=0$. Then there exist an elliptic curve $F$ and an integer $k\geq1$ such that $D^{'}\sim kF$. By Corollary 2.1 (i), we have $3k\leq H.D^{'}\leq H^2-2$. Since
$$\chi(\mathcal{O}_X(D^{'}))= h^0(\mathcal{O}_X(D^{'}))-k+1\geq \chi(\mathcal{O}_X(2H-D))-\frac{1}{3}(H^2-5),$$
\noindent we have the contradiction ${D^{'}}^2\geq \frac{1}{3}(H^2-8)>0$. Since ${D^{'}}^2>0$, we have ${D^{'}}^2\geq (2H-D)^2=H^2-6$. Since $H^2\geq10$, the Hodge index theorem implies that
$$(H^2-4)^2< H^2(H^2-6)\leq (H^2)({D^{'}}^2)\leq (H.D^{'})^2.$$
\noindent Hence, by the inequality (14), we have $H.D^{'}=H^2-3$ or $H^2-2$. 

Assume that $H.D^{'}=H^2-3$. If ${D^{'}}^2>H^2-6$, then $(H-D^{'})^2>0$. This implies that $(H-D^{'})^2\geq2$. Hence, by the Hodge index theorem, we have the contradiction 
$$9=(H.(H-D^{'}))^2\geq (H^2)(H-D^{'})^2\geq 20.$$
\noindent Hence, we have ${D^{'}}^2=H^2-6$. Since $H^2\geq10$, we have $h^1(\mathcal{O}_X(D^{'}))=0$. Since ${D^{'}}^2=(2H-D)^2$, by Lemma 3.2, we have the equality (12). 

Assume that $H.D^{'}=H^2-2$. Since $H.\Delta=1$, $\Delta$ is a $(-2)$-curve. Since $$(D^{'}+\Delta)^2=(2H-D)^2=H^2-6,$$
we have $2D^{'}.\Delta=H^2-4-{D^{'}}^2.$
Since $D^{'}.\Delta\geq0$, we have ${D^{'}}^2 = H^2-4$ or $H^2-6$. If ${D^{'}}^2=H^2-4$, we have $(H-D^{'})^2=0$ and $H.(H-D^{'})=2$. However, by Corollary 2.1 (i), this is a contradiction. Since ${D^{'}}^2=H^2-6$, we have $D^{'}.\Delta=1$. Since ${D^{'}}^2>0$, we may assume that $D^{'}$ is irreducible. Therefore, $|2H-D|$ contains the 1-connected divisor $D^{'}+\Delta$. Hence, we have the equality (12). $\hfill\square$

\begin{prop} Assume that $D^2\geq H^2$. Then $\mathcal{O}_X(D)$ is ACM and initialized with respect to $L$ if and only if $D^2=2H.D-H^2-4$, $|D-H|=\emptyset$, and $h^1(\mathcal{O}_X(2H-D))=0$.\end{prop}

{\it{Proof}}. First of all, we show that 
$$|H-D|=\emptyset.\leqno(15)$$
If $|H-D|\neq\emptyset$, by the ampleness of $H$, we have $H.(H-D)>0$, that is, $H.D<H^2$. By the Hodge index theorem, we have $H^2D^2\leq(H.D)^2$. Hence, we have $D^2<H^2$, which contradicts the hypotheses.

Assume that $\mathcal{O}_X(D)$ is ACM and initialized. By the assumption, we have $h^1(\mathcal{O}_X(2H-D))=0$ and $|D-H|=\emptyset$. Since $h^1(\mathcal{O}_X(H-D))=0$, by the equality (15), we have $\chi(\mathcal{O}_X(H-D))=0$, that is, $D^2=2H.D-H^2-4$.

Conversely, we assume that $D^2=2H.D-H^2-4$, $|D-H|=\emptyset$, and $h^1(\mathcal{O}_X(2H-D))=0$. By the equality (15), we have
$$h^1(\mathcal{O}_X(D-H))=0.$$
\noindent Since
$$\chi(\mathcal{O}_X(D-H))=0\text{ and }\chi(\mathcal{O}_X(2H-D))\geq0,$$
\noindent we have 
$$D^2\leq 2H^2-4.\leqno(16)$$
We show that $\mathcal{O}_X(D)$ is ACM. Since 
$$H.D=\frac{1}{2}(D^2+H^2+4)\leq\frac{3}{2}H^2,$$
\noindent by Lemma 3.1 and the previous vanishings, it is sufficient to show that 
$$h^1(\mathcal{O}_X(D))=0.$$
Since $D^2>0$, we show that $|D|$ is base point free. Assume that $|D|$ is not base point free. Let $\Delta$ be the fixed component of $|D|$, and let $D^{'}$ be the movable part of $|D|$. By the ampleness of $H$, we have
$$H.D^{'}<H.D\leq\frac{3}{2}H^2.$$
\noindent Assume that ${D^{'}}^2=0$. Then there exist an elliptic curve $F$ and a positive integer $k$ such that $D^{'}\sim kF$. By Proposition 2.4 (i), we have $H.D^{'}=kH.F\geq 3k.$
\noindent Since $k\leq\frac{1}{2}H^2-1$, we have
$$\chi(\mathcal{O}_X(D^{'}))=h^0(\mathcal{O}_X(D^{'}))-k+1\geq \chi(\mathcal{O}_X(D))-\frac{1}{2}H^2+2.$$
\noindent Hence, by the assumption that $D^2\geq H^2$, we have ${D^{'}}^2\geq D^2-H^2+4\geq 4.$
\noindent This is a contradiction. Hence, we have ${D^{'}}^2>0$. Since $h^1(\mathcal{O}_X(D^{'}))=0$, we have
$${D^{'}}^2\geq D^2\geq H^2.$$ 
\noindent By the Hodge index theorem, we have
$$(H^2)^2\leq (D^2)(H^2)\leq ({D^{'}}^2)(H^2)\leq (H.D^{'})^2.$$
\noindent Hence, we have $H.(D^{'}-H)\geq 0$. Since $H.\Delta>0$, by the assumption, we have 
$$(D^{'}-H)^2>(D-H)^2=-4.$$
\noindent This implies that $|D^{'}-H|\neq\emptyset$. However, this contradicts the assumption that $|D-H|=\emptyset$. Hence, $|D|$ is base point free. We have $h^1(\mathcal{O}_X(D))=0$. $\hfill\square$

\section{Example of ACM line bundles} Let the notations be as in Theorem 1.2. By the inequality (16), if $\mathcal{O}_X(D)$ is ACM and initialized with respect to $L$, then $D^2\leq 2H^2-4$. In particular, an ACM and initialized line bundle $\mathcal{O}_X(D)$ satisfying the equality $D^2=2H^2-4$ is called an Ulrich line bundle. In general, it is difficult to consider the problem whether such a line bundle exists or not, for a given polarization $L$. However, we can give an example of an initialized and ACM line bundle $\mathcal{O}_X(D)$ which is not an Ulrich line bundle and satisfies the condition as in Theorem 1.2 (ii) (d). First of all, we recall the following existence theorem.

\begin{thm}{\rm{(cf. [J-K], Proposition 4.2 and Lemma 4.3)}}. Let $d$ and $g$ be integers with $g\geq3$ and $3\leq d\leq \lfloor\displaystyle\frac{g+3}{2}\rfloor$. Then there exists a K3 surface $X$ with $\Pic(X)=\mathbb{Z}[C]\oplus\mathbb{Z}[F]$ such that $C$ is a smooth curve of genus $g$ and $F$ is an elliptic curve with $C.F=d$, where $[D]$ is the linearly equivalent class of a divisor $D$ on $X$. Moreover, $\Cliff(C)=d-2$.
\end{thm}

\noindent We note that, in Theorem 4.1, the gonality of $C$ is $d$ and it can be computed by a pencil on $C$ which is given by the restriction of the elliptic pencil $|F|$ on $X$ ([J-K], proof of Lemma 4.3 and Theorem 4.4). Moreover, we can see that $\Pic(X)$ contains a $(-2)$-vector if and only if $d|g$, by easy computation.

\begin{prop} Let $X$ be a K3 surface as in Theorem 4.1, and assume that $d|(g+1)$. If we let $m=\displaystyle\frac{g+1}{d}$, $L=\mathcal{O}_X(mC)$, and $D\sim (m+1)C-mF$, then $L$ is a very ample line bundle on $X$, and $\mathcal{O}_X(D)$ is ACM and initialized with respect to $L$.\end{prop}

{\it{Proof}}. Since $d|(g+1)$, we have $d\mid \hspace{-.67em}/g$. Since there is no $(-2)$-curve on $X$, $C$ is ample. Moreover, for any elliptic curve $E$, we have $C.E\geq d$. Indeed, if $E$ is an elliptic curve on $X$ which is not linearly equivalent to $F$, by easy computation, there exist two integers $s>0$ and $t$ such that $E\sim sC+tF$ and  $s(g-1)+td=0$. Since $d\neq g$, we have $\lfloor\displaystyle\frac{g+3}{2}\rfloor\leq g-1$. In fact, if $\lfloor\displaystyle\frac{g+3}{2}\rfloor> g-1$, then we have $g\leq3$. By the assumption, this implies that $d=g=3$. Hence, we have
$$E.C=s(g-1)\geq g-1\geq d.$$
Moreover, since $d\leq g-1$, we have $m\geq 2$. Since $L^2=m^2(2g-2)\geq16$, by Proposition 2.4, $L$ is very ample. Since $d\geq 3$, we have $g\geq 5$. Hence, we have $D^2-L^2=(2g-6)m-4>0$. Since $D.F>0$, we have $|D|\neq\emptyset$. Let $H\in |L|$. Since $g+1=md$, we have $(D-H)^2=-4$. Since there is no $(-2)$-curve on $X$, we have $|D-H|=\emptyset$. Obviously, $2H-D\sim (m-1)C+mF$ is base point free and big. Therefore, we have $h^1(\mathcal{O}_X(2H-D))=0$. By Theorem 1.2, $\mathcal{O}_X(D)$ is ACM and initialized with respect to $L$. $\hfill\square$

$\;$

\noindent {\bf{Acknowledgements}}. The author is partially supported by Grant-in-Aid for Scientific Research (16K05101), Japan Society for the Promotion Science.

\end{document}